\documentclass[12pt]{article}

\RequirePackage[italian,english]{babel}

\title{Characteristic curves for set-valued Hamilton-Jacobi equations}
\author{D.~Visetti\footnote{Free University of Bozen-Bolzano, 
  Faculty of Economics and Management, 
  \href{mailto:daniela.visetti@unibz.it}{daniela.visetti@unibz.it}.  
  The work was supported by the project \emph{Optimal Control Problems with 
  Set-valued Objective Function}, Free University of Bolzano-Bozen.}}

\RequirePackage[italian,english]{babel} 
\RequirePackage{graphicx}
\RequirePackage{amsmath}
\RequirePackage{amssymb}
\RequirePackage{amsthm}
\usepackage[colorlinks=true, linkcolor=blue, citecolor=blue, urlcolor=blue, filecolor=blue]{hyperref}
\usepackage{color}

\newcommand{\cl}{{\rm cl \,}}

\newcommand{\bs}{\backslash}
\renewcommand{\P}{\ensuremath{\mathcal{P}}}

\newcommand{\G}{\ensuremath{\mathcal{G}}}
\newcommand{\Ha}{\ensuremath{\mathcal{H}}}
\newcommand{\V}{\ensuremath{\mathcal{V}}}

\definecolor{color0}{gray}{.50}
\definecolor{color1}{rgb}{0,.2,.8}
\definecolor{color2}{rgb}{1,.2,0}
\definecolor{color3}{rgb}{.2,.7,.6}

\begin{document}
\selectlanguage{english}

\newcommand{\NN}{{\mathbb N}}
\newcommand{\RR}{{\mathbb R}}
\newcommand{\ZZ}{{\mathbb Z}}
\newcommand{\QQ}{{\mathbb Q}}
\newcommand{\CC}{{\mathbb C}}

\newtheorem{theorem}{Theorem}[section]
\newtheorem{proposition}[theorem]{Proposition}
\newtheorem{lemma}[theorem]{Lemma}
\newtheorem{corollary}[theorem]{Corollary}
\newtheorem{remark}[theorem]{Remark}
\newtheorem{example}[theorem]{Example}
\newtheorem{definition}[theorem]{Definition}

\maketitle

\begin{abstract}
  The method of characteristics is extended to set-valued Hamilton-Jacobi 
  equations.  This problems arises from a calculus of variations' problem with a 
  multicritiria Lagrangian function: through an embedding into a set-valued 
  framework, a set-valued Hamilton-Jacobi equation is derived, where 
  the Hamiltonian function is the Fenchel conjugate of the Lagrangian 
  function.  In this paper a method of characteristics is described and some 
  results are given for the Fenchel conjugate.
\end{abstract}


\bigskip


\section{Introduction}

The method of characteristics converts a nonlinear first-order partial differential 
equation to a system of ordinary differential equations, both with suitable 
boundary conditions.   In some circumstances, the solutions of the latter exist 
(at least locally in time) and give a solution of the first one.  This happens in 
particular for the Hamilton-Jacobi equation both for classical solutions and 
for viscosity solutions (see \cite{Lions}).  The characteristic curves are called 
generalized in the second case.

This method was already suggested by Cauchy and was mentioned in 
Cara\-th\'eodory's book \cite{Caratheodory}.   In order to have the general 
theory and some perspective in the research that has been done, it is possible 
to consider \cite{Evans}, \cite{Cannarsa_Sinestrari}, \cite{Subbotina}.

In this paper a Hamilton-Jacobi equation for a set-valued Hamiltonian function 
is presented.  The Hamiltonian function takes values in a complete lattice of sets 
that are invariant with respect to the sum of a closed and convex cone.  This 
problem was never studied before, to the knowledge of the author, and it is not 
only interesting per se, but also as a natural completion of a previous paper.   
In fact, in \cite{hv} the authors considered a calculus of variations' problem with 
a multicriteria Lagrangian function.  The embedding of it into a set-valued 
framework was fundamental to prove a Hopf-Lax formula for the value function.  
Subsequently, through the Fenchel conjugate of the Lagrangian function a 
Hamilton-Jacobi set-valued equation was derived.

Here, like in the classical theory, we start supposing that a smooth solution of 
the equation is known, in order to write the system of ordinary differential 
equations and the characteristic curves.  In this situation, a family of characteristic 
curves is obtained, parametrized by the elements of the dual cone.  A posteriori 
the ordinary differential equations do not depend on the solution and can be 
solved independently.   Under some assumptions it is so possible to write a 
solution of the set-valued Hamilton-Jacobi equation from the solutions of the 
ordinary differential equations.  The solutions are global under some hypotheses.

Considering the case of the Hamiltonian function as the Fenchel conjugate of a 
Lagrangian function, some results are obtained.  For set-valued optimization an 
infimizer is a set.  It is proved that an infimizer for the Hopf-Lax formula proved in 
\cite{hv} can be constituted only by minimizers of the scalarized problem with 
respect to an element of the dual cone.

Moreover, some properties regarding the derivatives of the Fenchel conjugate 
are extended to the set-valued framework.  Finally, it is proved that the characteristic curves coincide with elements of the infimizer.

This study opens a variety of new questions.  First of all, it is still not clear what 
can be an extension of the concept of viscosity solution for a set-valued 
Hamilton-Jacobi equation.  Then also the problem of the generalized characteristic 
curves should be addressed.  Recently, in \cite{nh} 
the author considers a generalized form of characteristics.  But also in the 
previously cited texts \cite{Cannarsa_Sinestrari} and \cite{Subbotina}, there 
are several approaches that have not yet been studied in the set-valued case.

\section{Preliminaries}

The Minkowski sum of two non-empty sets $A, B \subseteq \RR^d$ is $A + B = \{a + b \mid a \in A, \, b \in B\}$. It is  
extended to the whole power set $\mathcal{P}(\RR^d)$ by
$$
\emptyset + A = A+ \emptyset= \emptyset.
$$
We also use $A \oplus B := \cl(A+B)$, the `closed sum' of two sets.

A set $C \subseteq\RR^d$ is a cone if $s C \subseteq C$ for all $s > 0$, and it is a 
convex cone if additionally $C + C \subseteq C$.    Let $C$ be a closed and 
convex cone in $\RR^d$ different from the empty set and the 
whole $\RR^d$.  The dual of a cone $C$ is defined as
$$
C^+=\{\zeta\in\RR^d \mid \zeta\cdot z\geq 0\}.
$$
If there is an element $\hat z \in C$ such that $\zeta\cdot\hat z > 0$ for all $\zeta\in 
C^+\bs\{0\}$ (in particular, if $\mbox{int}\; C \neq \emptyset$), then the set
\begin{equation}
\label{base}
B^+(\hat z) = \{\zeta \in C^+ \mid \zeta\cdot\hat z = 1\}
\end{equation}
is a (closed and convex) base of $C^+$, i.e., for each element $\xi\in C^+\bs\{0\}$ 
there are unique $\zeta\in B^+(\hat z)$ and $s > 0$ such that $\xi = s\zeta$.

We consider the following subset of the power set $\mathcal{P}(\RR^d)$ 
(see for instance \cite{SetOptSurvey}):
$$
\mathcal{G}(\RR^d,C)=\{A\in\mathcal{P}(\RR^d)\ \vert\ A = \mathrm{cl\, co}\,
  (A+C)\}
$$
where cl and co are the closure and the convex hull, respectively.

The pair $(\mathcal{G}(\RR^d,C),\supseteq)$ is a complete lattice. If $\mathcal{A}
\subseteq\mathcal{G}(\RR^d,C)$, then the infimum and the supremum of 
$\mathcal{A}$ are given by
$$
\inf\mathcal{A}=\mathrm{cl\, co}\bigcup_{A\in\mathcal{A}}A\qquad
\sup\mathcal{A}=\bigcap_{A\in\mathcal{A}}A.
$$
An element $A_0\in\mathcal{A}$ is called minimal for $\mathcal{A}$ if 
$$
A\in\mathcal{A},\ A\supseteq A_0\ \Longrightarrow\ A=A_0\, .
$$

Let $\zeta\in C^+\bs\{0\}$ and let
$$
H^+(\zeta)=\{z \in \RR^d \mid \zeta \cdot z\geq 0 \}
$$
where $\zeta \cdot z$ denotes the usual scalar product.
For two sets $A,B\in\mathcal{P}(\RR^d,C)$, the set
$$
A-_{\zeta}B=\{ z\in\RR^d \mid z+B \subseteq A \oplus H^+(\zeta)\}
$$
is called the $\zeta$-difference of $A$ and $B$.  It is either $\emptyset$, $\RR^d$ or a closed (shifted) half-space.  It is possible to see that
$$
A -_{\zeta} B =\{z\in \RR^d \mid \zeta \cdot z+\inf_{b\in B} \zeta \cdot b \geq
    \inf_{a\in A} \zeta \cdot a\}
$$

Let $\{ A_n\}_{n\in\NN}$ be a sequence of sets in $\mathcal{G}(\RR^d,C)$, we 
denote by $\lim_{n\to\infty}A_n$ the following set:
$$
\lim_{n\to\infty}A_n=\left\{ z\in\RR^d\mid\forall n\in\NN, \exists
  z_n\in A_n, \lim_{n\to\infty}z_n=z\right\}.
$$
This definition of limit coincides with the upper limit of Painlev\'e-Kuratowski 
($\mathrm{Liminf}_{n\to\infty}A_n=\left\{ z\in Z\ \vert\ \lim_{n\to\infty}d(z,A_n)
=0\right\}$, see \cite{Aubin_Frank}).

Let $\{ A_s\}_{s\in S}$ with $S\subseteq\RR$ be a family of sets in $\mathcal{G}
(\RR^d,C)$ and $\bar s \in \RR$. We denote by $\lim_{s \to \bar s}A_s$ the set 
which satisfies that for any sequence $\{s_n\}_{n\in\NN}\subseteq S$ with $s_n 
\to \bar s$ one has
$$
\lim_{s \to \bar s}A_s = \lim_{n\to\infty}A_{s_n}.
$$

Let $f$ be a function $f:\RR^n\to\mathcal{G}(\RR^d,C)$.  The graph of $f$ is 
$$
\mbox{graph}\ f=\{(x,z)\in \RR^n\times\RR^d\mid z\in f(x)\}\, .
$$
The domain of $f$ is
$$
\mbox{dom}\ f=\{x\in\RR^n\mid f(x)\neq\emptyset\}
$$
The function is convex if and only if 
$\mbox{graph}\ f$ is a convex subset of $\RR^n\times\RR^d$.  This is equivalent 
to the following condition: for any $\lambda\in(0,1)$, $x_1,x_2\in\RR^n$
$$
f(\lambda x_1+(1-\lambda)x_2)\supseteq \lambda f(x_1)+(1-\lambda)f(x_2)\, .
$$

Let $(X,\mathcal{A},\mu)$ be a measure space and let $f$ be a set-valued map 
from $X$ into the closed nonempty subsets of $\RR^d$.

The set of the integrable selections of $f$ is:
$$
\mathcal{F}=\{\varphi\in L^1(X,\RR^d)\ |\ \varphi(x)\in f(x) \mbox{ a.e. in }X\}
$$

The Aumann integral of $f$ on $\RR^n$ is the set of integrals of the integrable 
selections of $f$:
$$
\int_X f\; d\mu = \left\{ \int_{\RR^n} \varphi\; d\mu\ |\ \varphi\in
  \mathcal{F}\right\}
$$

Let $X$ be a non-empty set, $f:X\to\mathcal{G}(\RR^d,C)$ a function and $f[X]=
\{ f(x)\ |\ x\in X\}$.

A family of extended real-valued functions $\varphi_{f, \zeta} \colon X \to\overline
\RR$ with $\zeta\in C^+$ is defined by
\begin{equation}
\label{EqScalarization}
\varphi_{f,\zeta}(x) = \inf_{z \in f(x)} \zeta\cdot z.
\end{equation}
A point $\bar x \in X$ is called a $\zeta$-minimizer of $f$ if 
$$
\forall x \in X \colon \varphi_{f,\zeta}(\bar x) \leq \varphi_{f,\zeta}(x).
$$

\begin{itemize}
\item[(a)] A set $M\subset X$ is called an \emph{infimizer for $f$} if
  $$
  \inf f[M]=\inf f[X]\, .
  $$
\item[(b)] An element $x_0\in X$ is called a \emph{minimizer for $f$} if
  $f(x_0)$ is minimal for $f[X]$.
\item[(c)] A set $M\subset X$ is called a \emph{solution} of the problem 
  \emph{minimize $f(x)$ subject to $x\in X$} if $M$ is an infimizer for $f$ 
  and each $x_0\in M$ is a minimizer for $f$.  It is called a \emph{full solution} 
  if the set $f[M]$ includes all minimal elements of $f[X]$.
\item[(d)] A set $M \subseteq X$ is called a \emph{scalarization solution} of 
  the problem \emph{minimize $f(x)$ subject to $x\in X$} if it is an infimizer 
  and only includes $\zeta$-minimizers.
\end{itemize}

The solution concept in (d) has been considered first in \cite{AC}.

Let $\eta\in\RR^n$ and $\zeta \in C^+$ be given. We recall the definition of the 
function $S_{(\eta,\zeta)}:\RR^n\to\mathcal{G}(\RR^d,C)$:
$$
S_{(\eta,\zeta)}(x)=\{z\in\RR^d\ \vert\ \zeta\cdot z\geq\eta\cdot x\}\, .
$$
Such a function is additive and positively homogeneous, i.e., for all $x\in\RR^n$, 
$\lambda>0$
$$
S_{(\eta,\zeta)}(\lambda x)=\lambda S_{(\eta,\zeta)}(x)
$$
and for all $x_1,\, x_2\in\RR^n$
$$
S_{(\eta,\zeta)}(x_1+x_2)=S_{(\eta,\zeta)}(x_1)+S_{(\eta,\zeta)}(x_2)\, .
$$
Let $\hat z\in\RR^d$ be such that $\zeta\cdot\hat z=1$.  Then for any 
$x\in\RR^n$
\begin{equation}
\label{repr_S}
S_{(\eta,\zeta)}(x)=(\eta\cdot x)\hat z+H^+(\zeta)
\end{equation}
(see \cite{SetOptSurvey}).

A derivative concept that is most compatible with the complete lattice approach 
is as follows.   The derivatives of a function $f:\RR\times\RR^n\to\P(\RR^d,C)$ 
with respect to an element $\zeta$ in the dual cone, if they exist, will be defined 
in the following way:
\begin{equation}
\label{set_der}
\begin{aligned}
D_{\zeta,t}f(t,x) &= \lim_{h\to 0}\frac1h[f(t+h,x)-_\zeta f(t,x)], \\
D_{\zeta,x}f(t,x)(q) &= \lim_{h\to 0} \frac1h[f(t,x+hq)-_\zeta f(t,x)].
\end{aligned}
\end{equation}
See \cite{AC} and \cite{hv} for a motivation and many features including a 
discussion of the `improper' function values $\RR^d$ and $\emptyset$.

The Fenchel conjugate of the function $f:\RR^n\to\mathcal{P}(\RR^d,C)$ is 
defined as the function
\begin{equation}
\label{Fenchel}
\begin{array}{cccc}
f^*: & \RR^n\times C^+\backslash\{0\} & \to & \mathcal{G}(\RR^d,C) \\
       & (\eta,\zeta) & \mapsto & \sup_{x\in\RR^n} S_{(\eta,\zeta)}(x)-_\zeta f(x)
\end{array}
\end{equation}
%

\section{Hamilton-Jacobi equation and characteristic curves}

Let us consider
$$
0<T< +\infty,\quad Q_T=[0,T]\times\RR^n, \quad (t,x)\in\RR\times\RR^n,
$$
a set-valued Hamiltonian function $\Ha:\RR^n\times C^+\bs\{0\}\to\G(\RR^d,C)$
and $U_0:\RR^n\to\RR^d$ a function of class $C^2$.     The form of this 
Hamiltonian is due to the fact that sometimes it 
is the Fenchel conjugate \eqref{Fenchel} of a Lagrangian function (see \cite{hv}).

We suppose that for any $\zeta\in C^+\bs\{0\}$, the function
$\Ha_\zeta:\RR^n\to\RR$, defined as
$$
\Ha_\zeta(p)=\inf_{z\in\Ha(p,\zeta)}\zeta\cdot z,
$$
be of class $C^2$.   

Another particular case can be if there exists $\Ha_0:\RR^n\to\RR^d$ of class 
$C^2$ and if $\Ha$ is in some sense the inf-extension of $\Ha_0$, $\Ha(p,\zeta)=
\Ha_0(p)+C$ (see \cite{SetOptSurvey}), then $\Ha_\zeta(p)$ coincide 
with $\Ha_{0,\zeta}(p)=\Ha_0(p)\cdot\zeta$ and are automatically of class $C^2$.

Given $\zeta\in B^+(\hat z)$ (see \eqref{base}), let $U:Q_T\to
\G(\RR^d,C)$ be such that $U(t,x)+H^+(\zeta)$ can 
be written as
\begin{equation}
\label{u_zeta}
U(t,x)+H^+(\zeta)=u_\zeta(t,x)\hat z+H^+(\zeta)\, ,
\end{equation}
with $u_\zeta:Q_T\to\RR$.  
Property (\ref{u_zeta}) gives (see \eqref{set_der}):
$$
\begin{aligned}
D_{\zeta,t}U(t,x) &= S_{\left( \frac{\partial u_\zeta}{\partial t}(t,x),\zeta\right)}(1)\, , \\
D_{\zeta,x}U(t,x)(q) &= S_{\left(Du_\zeta(t,x),\zeta\right)}(q)\, .
\end{aligned}
$$

For $\zeta\in B^+(\hat z)$, we suppose that the Hamilton-Jacobi equation
\begin{equation}
\label{Ham_Jac_eq_zeta}
\left\{ \begin{array}{l} 
  D_{\zeta,t}U(t,x)+\Ha\left(Du_{\zeta}(t,x),\zeta\right)=
    H^+(\zeta) \\
  U(0,x)=U_0(x)+C
  \end{array} \right.
\end{equation}
admits a solution $U(t,x)$ on $Q_T$ with the property (\ref{u_zeta}) and where 
$u_\zeta$ is of class $C^2$.  If $U(t,x)$ is a solution of \eqref{Ham_Jac_eq_zeta} for every $\zeta\in B^+(\hat z)$, then it is also a solution of 
\begin{equation}
\label{Ham_Jac_eq}
\left\{ \begin{array}{l} 
  \displaystyle\sup_{\zeta\in B^+(\hat z)} \left[D_{\zeta,t}U(t,x)+
    \Ha\left(Du_{\zeta}(t,x),\zeta\right)\right]=C \\
  U(0,x)=U_0(x)+C
  \end{array} \right.
\end{equation}

We say that $U(t,x)$, satisfying property (\ref{u_zeta}), is of class $C^2$ if all 
the $u_\zeta$ are $C^2$ for $\zeta\in B^+(\hat z)$.

The first equation in (\ref{Ham_Jac_eq_zeta}) gives
$$
\inf\left\{\zeta\cdot z\mid z\in\left[ D_{\zeta,t}U(t,x)+\Ha
  \left(Du_{\zeta}(t,x),\zeta\right) \right]\right\} = \inf\{\zeta\cdot z\mid 
  z\in H^+(\zeta)\}
$$
that can be written
$$
\frac{\partial u_\zeta}{\partial t}(t,x)+\mathcal{H}_\zeta\left(Du_\zeta
  (t,x)\right)=0\, .
$$
So $u_\zeta$ is a solution of
\begin{equation}
\label{Ham_Jac_eq_scal}
\left\{ \begin{array}{l}
\frac{\partial u_\zeta}{\partial t}(t,x)+\mathcal{H}_\zeta\left(Du_\zeta
  (t,x)\right)=0 \\
  u_\zeta(0,x)=U_{0,\zeta}(x)
  \end{array} \right.
\end{equation}
where $U_{0,\zeta}(x)=\zeta\cdot U_0(x)$ is a real-valued function.

For fixed $x\in\RR^n$, we denote by $X_\zeta(t,x)$ the solution of the ordinary 
differential equation
\begin{equation}
\label{X_zeta}
\dot X_\zeta=D\mathcal{H}_\zeta\left(Du_\zeta(t,X_\zeta)\right), 
  \qquad X_\zeta(0,x)=x\, .
\end{equation}
Such a solution is defined on an interval $[0,T_{\zeta,x})$.  The curve
$(t,X_\zeta(t,x))$ is the characteristic curve associated to $U$ with respect to 
$\zeta$.

We define now
\begin{equation}
\label{V_zeta,P_zeta}
V_\zeta(t,x)=u_\zeta(t,X_\zeta(t,x))\, , \qquad 
P_\zeta(t,x)=Du_\zeta(t,X_\zeta(t,x))\, .
\end{equation}
Using equation (\ref{Ham_Jac_eq_scal}), we obtain
$$
\begin{aligned}
\dot V_\zeta &= \frac{\partial u_\zeta}{\partial t}(t,X_\zeta)+D
  u_\zeta(t,X_\zeta)\cdot\dot X_\zeta = - \mathcal{H}_\zeta\left( P_\zeta
  \right) + D\mathcal{H}_\zeta\left( P_\zeta \right)\cdot P_\zeta\, , \\
\dot P_\zeta &= \frac{\partial Du_\zeta}{\partial t}(t,X_\zeta)+
  D^2u_\zeta(t,X_\zeta)\dot X_\zeta \\
  &= D \left( \frac{\partial 
  u_\zeta}{\partial t}(t,X_\zeta) + \mathcal{H}_\zeta\left(Du_\zeta
  (t,X_\zeta)\right)\right)=0\, .
\end{aligned}
$$
As a consequence, $P_\zeta$ is constant in time:
$$
P_\zeta(t,x) \equiv DU_{0,\zeta}(x)\, .
$$
Then also $\dot X_\zeta$ is constant in time: $\dot X_\zeta=D\mathcal{H}_\zeta\left(P_\zeta\right)=D\mathcal{H}_\zeta\left(DU_{0,\zeta}(x)\right)$.   The 
solutions of the ODEs are then:
\begin{equation}
\label{XVP}
\left\{ \begin{array}{l}
  X_\zeta(t,x) = x+tD\mathcal{H}_\zeta\left(DU_{0,\zeta}(x)\right) \\
  V_\zeta(t,x) = U_{0,\zeta}(x)+t\left( - \mathcal{H}_\zeta\left(DU_{0,\zeta}(x)
    \right)+ D\mathcal{H}_\zeta\left(DU_{0,\zeta}(x)\right)\cdot
    DU_{0,\zeta}(x)\right) \\
  P_\zeta(t,x) = DU_{0,\zeta}(x)
\end{array} \right.
\end{equation}

The next step is to consider the system of ODEs
\begin{equation}
\label{ODE}
\left\{ \begin{array}{l}
\dot X_\zeta=D\mathcal{H}_\zeta\left(P_\zeta\right) \\
\dot P_\zeta=0 \\
\dot V_\zeta= - \mathcal{H}_\zeta\left( P_\zeta
  \right) + D\mathcal{H}_\zeta\left( P_\zeta \right)\cdot P_\zeta
\end{array} \right.
\end{equation}
in order to build a solution of the Hamilton-Jacobi equation.
The classical result is the following \emph{local existence theorem} (see for example 
\cite{Cannarsa_Sinestrari}):

\bigskip
\emph{
  \textit{Let $U_0$ be in $C^2(\RR^n;\RR^d)$ and  
  $DU_{0,\zeta}$, $D^2 U_{0,\zeta}$ be bounded for any $\zeta\in C^+
  \bs\{0\}$, $\|\zeta\|=1$. Let $\Ha_\zeta$ be of class $C^2$ for any $\zeta
  \in C^+\bs\{0\}$, $\|\zeta\|=1$. Denoting
  \begin{equation}
  \label{T*zeta}
  T^*_\zeta=\sup\{ t>0\ |\ I+tD^2\mathcal{H}_\zeta\left(DU_{0,\zeta}(x)
    \right)D^2U_{0,\zeta}(x) \mbox{ is invertible }\forall x\in\RR^n\}\, ,
  \end{equation}
  problem (\ref{Ham_Jac_eq_scal}) has a unique solution $u_\zeta\in C^2([0,
  T^*_\zeta)\times\RR^n)$.}
}

\bigskip

By the previous hypotheses, for any $T<T^*_\zeta$, there exists $Z_\zeta:[0,T]
\times\RR^n\to\RR^n$ of class $C^1$ such that
$$
X_\zeta(t,Z_\zeta(t,x))=x
$$
for any $(t,x)\in[0,T]\times\RR^n$.  The solution of the previous theorem is 
given by
\begin{equation}
\label{uzeta}
u_\zeta(t,x)=V_\zeta(t,Z_\zeta(t,x))\, ,\qquad \mbox{for any }(t,x)\in[0,T]
  \times\RR^n\, .
\end{equation}

We want now to find set-valued functions $U(t,x)$ that are solutions of the 
set-valued Hamilton-Jacobi equations.

\begin{theorem}
\label{main_trm}
  Let $U_0$ be in $C^2(\RR^n;\RR^d)$ and 
  $DU_{0,\zeta}$, $D^2 U_{0,\zeta}$ be bounded for any $\zeta\in B^+(\hat z)$.  
  Let $\Ha_\zeta$ be in $C^2(\RR^n)$ for any $\zeta\in B^+(\hat z)$.  Denoting
  \begin{equation}
  \label{T*}
  T^*=\inf_{\zeta\in B^+(\hat z)} T^*_\zeta
  \end{equation}
  where $T^*_\zeta$ is as in (\ref{T*zeta}), suppose that $T^*>0$. For 
  $T<T^*$ the map $U_\zeta:[0,T]\times\RR^n\to\mathcal{G}(\RR^n,C)$ 
  defined as
  $$
  U_\zeta(t,x) =u_\zeta(t,x)\hat z+H^+(\zeta)
  $$
  is a solution of
  \begin{equation}
  \label{Ham_Jac_eq_zeta_H}
  \left\{ \begin{array}{l} 
  D_{\zeta,t}U(t,x)+\Ha\left(Du_{\zeta}(t,x),\zeta\right)=
    H^+(\zeta) \\
  U(0,x)=U_0(x)+H^+(\zeta).
  \end{array} \right.
  \end{equation}

  Moreover, let $U$ be the map $U:[0,T]\times\RR^n\to\mathcal{G}(\RR^n,C)$ 
  defined as
  \begin{equation}
  \label{U}
  U(t,x)=\bigcap_{\zeta\in B^+(\hat z)} U_\zeta(t,x).
  \end{equation}
  If for any $\zeta\in B^+(\hat z)$ and $(t,x)\in[0,T]\times\RR^n$, there holds
  \begin{equation}
  \label{hyp_U}
  \inf_{z\in U(t,x)}\zeta\cdot z=\inf_{z\in U_\zeta(t,x)}\zeta\cdot z=u_\zeta(t,x),
  \end{equation}
  then $U$ is a solution of (\ref{Ham_Jac_eq}).
\end{theorem}

\begin{proof}
  It is easy to see that $U_\zeta(t,x)$ is a solution of  \eqref{Ham_Jac_eq_zeta_H}.
  
  The map $U(t,x)$ has the property (\ref{u_zeta}) thanks to hypothesis 
  \eqref{hyp_U}.

  The function $U$ satisfies the first equation in (\ref{Ham_Jac_eq}), because 
  $D_{\zeta, t}U(t,x)=D_{\zeta, t}U_\zeta(t,x)$. In fact, for any $A,B\in\G(\RR^d,C)$
  $$
  A-_\zeta B=(A+H^+(\zeta))-_\zeta(B+H^+(\zeta))
  $$
  and so also for the derivative
  $$
  D_{\zeta, t}U(t,x)=D_{\zeta, t}[U(t,x)+H^+(\zeta)]=D_{\zeta,t} U_\zeta(t,x).
  $$

  For the initial condition, we have that
  $$
  U(0,x)=\bigcap_{\zeta\in B^+(\hat z)}  U_\zeta(0,x)=
    \bigcap_{\zeta\in B^+(\hat z)} [U_0(x)+H^+(\zeta)]=U_0(x)+C\, .
  $$
\end{proof}

In the following proposition the solution is written as a set-valued version of the 
characteristic method.    In order 
to do that, we define a set-valued correspondent to $V_\zeta(t,x)$ for $\zeta\in 
B^+(\hat z)$.  More precisely, we denote $\V_\zeta:[0,T]\times\RR^n\to
\G(\RR^d,C)$ the map
\begin{equation}
\label{V_zeta_set}
\begin{aligned}
\V_\zeta(t,x) = &U_0(x)+t\left[ (H^+(\zeta)-_\zeta \Ha\left(DU_{0,\zeta}(x),\zeta
      \right))\right. \\
    &+ \left. D_{\zeta,p}\Ha\left(DU_{0,\zeta}(x),\zeta\right)(DU_{0,\zeta}(x))\right],
\end{aligned}
\end{equation}
where the derivative $D_\zeta$ denotes the derivative of the set-valued function 
$\Ha$.

\begin{proposition}
  Let $U_0$ be in $C^2(\RR^n;\RR^d)$ and 
  $DU_{0,\zeta}$, $D^2 U_{0,\zeta}$ be bounded for any $\zeta\in B^+(\hat z)$.  
  Let $\Ha_\zeta$ be in $C^2(\RR^n)$ for any $\zeta\in B^+(\hat z)$.  If 
  \eqref{hyp_U} holds, for $T^*>0$ as in Theorem \ref{main_trm} the solution 
  $U:[0,T]\times\RR^n\to\G(\RR^n,C)$ with $T<T^*$ defined in \eqref{U} can 
  be written
  $$
  U(t,x)=\bigcap_{\zeta\in B^+(\hat z)}\V_\zeta(t,Z_\zeta(t,x)).
  $$
\end{proposition}

\begin{proof}
  Both $U_\zeta(t,x)$ and $\V_\zeta(t,Z_\zeta(t,x))$ are half-spaces in the positive 
  direction of $\zeta$.  We have that
  $$
  \begin{aligned}
  \inf_{z\in H^+(\zeta)-_\zeta \Ha\left(DU_{0,\zeta}(x),\zeta\right)} \zeta\cdot z 
    &= -\Ha_\zeta\left(DU_{0,\zeta}(x)\right), \\
  \inf_{z\in D_{\zeta,p}\Ha\left(DU_{0,\zeta}(x),\zeta\right)(DU_{0,\zeta}(x))} 
    \zeta\cdot z &= D\Ha_\zeta\left(DU_{0,\zeta}(x)\right)\cdot DU_{0,\zeta}(x).
  \end{aligned}
  $$
  It is then immediate that $U_\zeta(t,x)=\V_\zeta(t,Z_\zeta(t,x))$.
\end{proof}

Next corollary presents a special case in which condition \eqref{hyp_U} in 
Theorem~\ref{main_trm} is satisfied.

\begin{corollary}
  Let $U_0$ be in $C^2(\RR^n;\RR^d)$ and 
  $DU_{0,\zeta}$, $D^2 U_{0,\zeta}$ be bounded for any $\zeta\in B^+(\hat z)$.  
  Let $\Ha(\cdot,\zeta)=\Ha_0(\cdot)+C$ for any $\zeta\in B^+(\hat z)$ and $\Ha_0$ 
  be in $C^2(\RR^n)$.  For $T^*>0$ as in Theorem \ref{main_trm} and $T<T^*$, 
  if for any $\zeta,\xi\in B^+(\hat z)$
  \begin{equation}
  \label{hyp_U_2}
  \begin{aligned}
  U_{0,\zeta} &(Z_\xi)+t\left[-\Ha_{0,\zeta}(DU_{0,\xi}(Z_\xi))+D\Ha_0(DU_{0,\xi}
      (Z_\xi))DU_{0,\xi}(Z_\xi)\right] \\
    &\geq U_{0,\zeta} (Z_\zeta)+t\left[-\Ha_{0,\zeta}(DU_{0,\zeta}(Z_\zeta))+
      D\Ha_0(DU_{0,\zeta}(Z_\zeta))DU_{0,\zeta}(Z_\zeta)\right],
  \end{aligned}
  \end{equation}
  where $Z_\zeta=Z_\zeta(t,x)$, $Z_\xi=Z_\xi(t,x)$ and $D\Ha_0$ denotes the 
  Jacobian matrix, then the function $U(t,x)$ 
  defined in \eqref{U} is a solution of \eqref{Ham_Jac_eq}.
\end{corollary}

\begin{proof}
  We observe that
  $$
  \V_\zeta(t,x)=U_0(x)+t\left[-\Ha_0(DU_{0,\zeta}(x))+D\Ha_0(DU_{0,\zeta}(x))
    DU_{0,\zeta}(x)\right] + H^+(\zeta).
  $$
  Equation \eqref{hyp_U_2} implies that the vectors
  $$
  U_0(Z_\zeta)+t\left[-\Ha_0(DU_{0,\zeta}(Z_\zeta))+D\Ha_0(DU_{0,\zeta}(Z_\zeta))
    DU_{0,\zeta}(Z_\zeta)\right]\in U(t,x)
  $$
  and that hypothesis \eqref{hyp_U} is satisfied.
\end{proof}

\begin{remark}
  In the hypotheses of Theorem \ref{main_trm}, we recall the following 
  results linked to the existence of the characterisics of the scalarized problems.
  \begin{enumerate}
  \item Set
    $$
    \begin{aligned}
    M_0 &= \sup_{\zeta\in B^+(\hat z)}
      \sup_{x\in\RR^n}\| DU_{0,\zeta}(x)\|, \\
    M_1 &= \sup_{\zeta\in B^+(\hat z)}
      \sup_{x\in\RR^n}\| D^2U_{0,\zeta}(x)\|, \\
    M_2 &= \sup_{\zeta\in B^+(\hat z)}
      \sup_{x\in\RR^n,\ \|x\|\leq M_0}\| D^2\mathcal{H}_\zeta(x)\|.
    \end{aligned}
    $$
    Then problem \eqref{Ham_Jac_eq} has a $C^2$ solution at least for the time 
    $t\in\left[0,\frac{1}{M_1M_2}\right)$.
  \item If $U_0$ and $\Ha_\zeta$ are convex, then problem \eqref{Ham_Jac_eq} 
    has a $C^2$ solution for all positive times.
  \item If $U_0(x)=Ax+b$ with $A$ a matrix of dimension $d$ times $n$ and 
    $b\in\RR^d$, then problem \eqref{Ham_Jac_eq} has a $C^2$ solution for all 
    positive times.
  \end{enumerate}
\end{remark}

\begin{example}
  We choose $\RR^n=\RR^d=\RR^2$, $C=\RR^2_+$ and the following Hamiltonian 
  and initial condition
  $$
  \Ha_0(p)=\left(\begin{array}{cc}
    \frac12\| p\|^2 \\ \frac14\| p\|^4 \end{array}\right), \qquad 
  \Ha(p,\zeta)= \Ha_0(p)+C, \qquad
  U_0(x)=\left(\begin{array}{cc} 1 & 0 \\ 0 & -1 \end{array}\right) x.
  $$
  If we consider $\hat z={1\choose 1}$, $B^+(\hat z)$ is given by all $\zeta=
  {\zeta_1\choose1-\zeta_1}$ for $0\leq\zeta_1\leq1$.  Since
  $$
  \begin{aligned}
  \Ha_\zeta(p) &= \frac{\zeta_1}{2}\| p\|^2 +\frac{1-\zeta_1}{4}\| p\|^4, \\
  D\Ha_\zeta(p) &= \zeta_1p+(1-\zeta_1)\| p\|^2p, \\
  U_{0,\zeta}(x) &= \zeta\cdot Ax, \\
  DU_{0,\zeta}(x) &= A^T\zeta,
  \end{aligned}
  $$
  where $A^T$ denotes the transpose of the matrix $A$.  
  Solving the system of ODEs \eqref{ODE}, we obtain
  $$
  \left\{ \begin{array}{l}
  X_\zeta(t,x) = x+t\left[\zeta_1+(1-\zeta_1)\| A^T\zeta\|^2\right]A^T\zeta \\
  V_\zeta(t,x) = \zeta\cdot Ax+t\left[\frac12\zeta_1\| A^T\zeta\|^2+
    \frac34(1-\zeta_1)\| A^T\zeta\|^4\right] \\
  P_\zeta(t,x) = A^T\zeta
  \end{array} \right.
  $$
  We want to calculate $Z_\zeta(t,x)$ such that
  $$
  X_\zeta(t,Z_\zeta(t,x))=x,
  $$
  which gives
  $$
  Z_\zeta(t,x)=x-t\left[\zeta_1+(1-\zeta_1)\| A^T\zeta\|^2\right]A^T\zeta.
  $$
  Now we can calculate
  $$
  U_\zeta(t,x)=Ax-t[\zeta_1+(1-\zeta_1)\| A^T\zeta\|^2]AA^T\zeta
    +t{\frac12\| A^T\zeta\|^2\choose\frac34\| A^T\zeta\|^4}+H^+(\zeta).
  $$
  In particular, for $x_0={1\choose 2}$ and $t=1$ the curve $\gamma_\zeta(1,x_0)=
  Ax_0-[\zeta_1+(1-\zeta_1)\| A^T\zeta\|^2]AA^T\zeta+{\frac12\| A^T\zeta\|^2
  \choose\frac34\| A^T\zeta\|^4}$ is plotted in the following figure:\\
  \begin{center}
  \includegraphics[width=10cm]{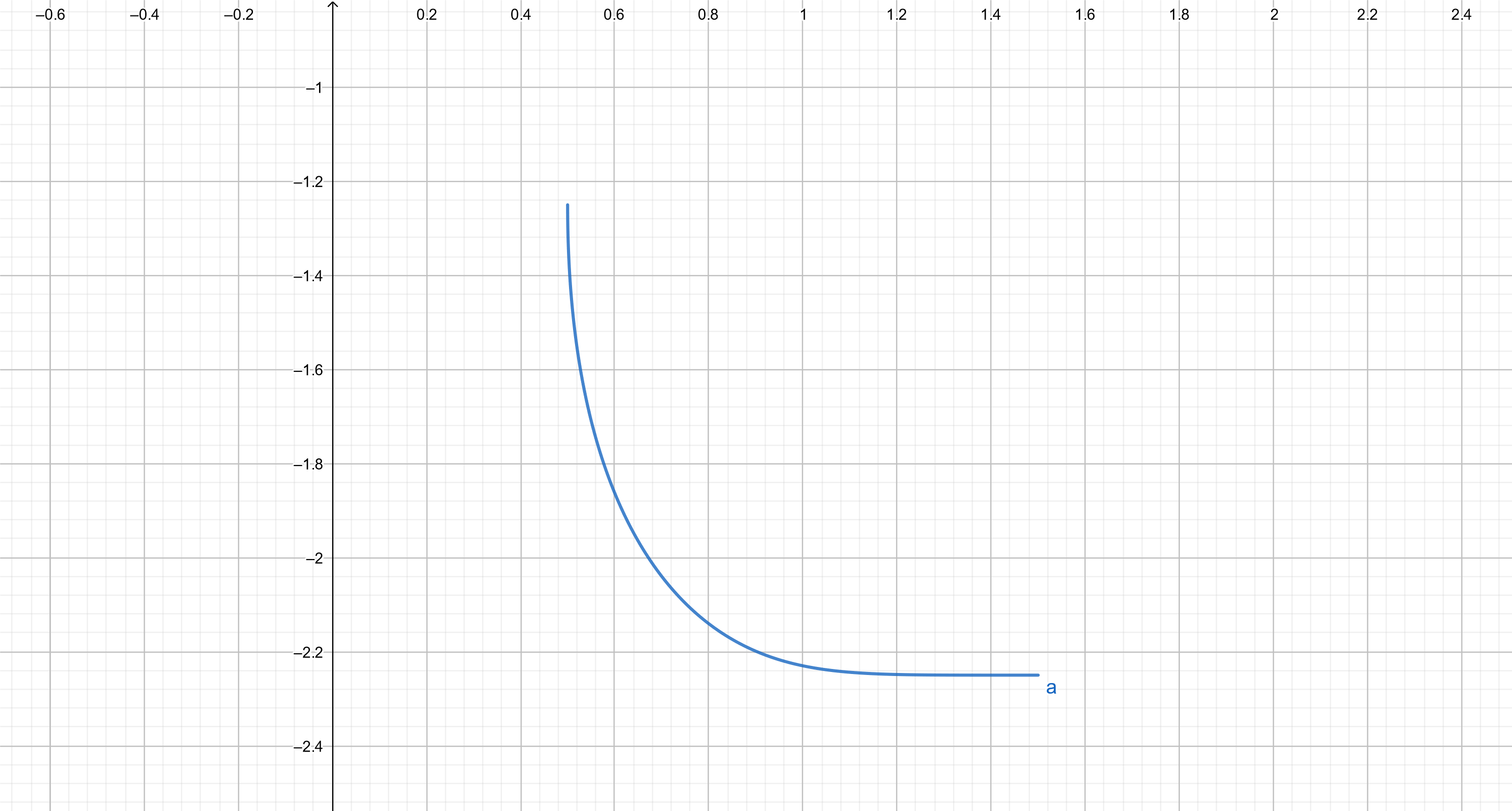}
  \end{center}
  while in the next two figures there are some half-spaces (corresponding to $\zeta=
  {1\choose0}$, ${3/(3+\sqrt3)\choose\sqrt3/(3+\sqrt3)}$, 
  ${1/2\choose1/2}$, ${1/(1+\sqrt3)\choose\sqrt{3}/(1+\sqrt3)}$, 
  ${0\choose1}$) and their intersection, which is an approximation of the 
  corresponding solution $U(1,x_0)$.\\
  \begin{center}
  \includegraphics[width=6cm]{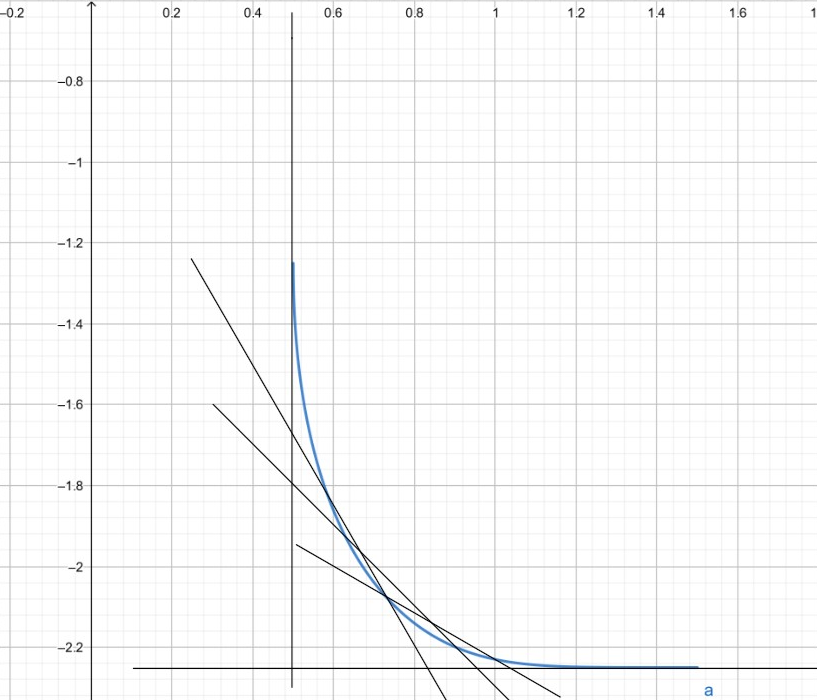}\quad
  \includegraphics[width=6cm]{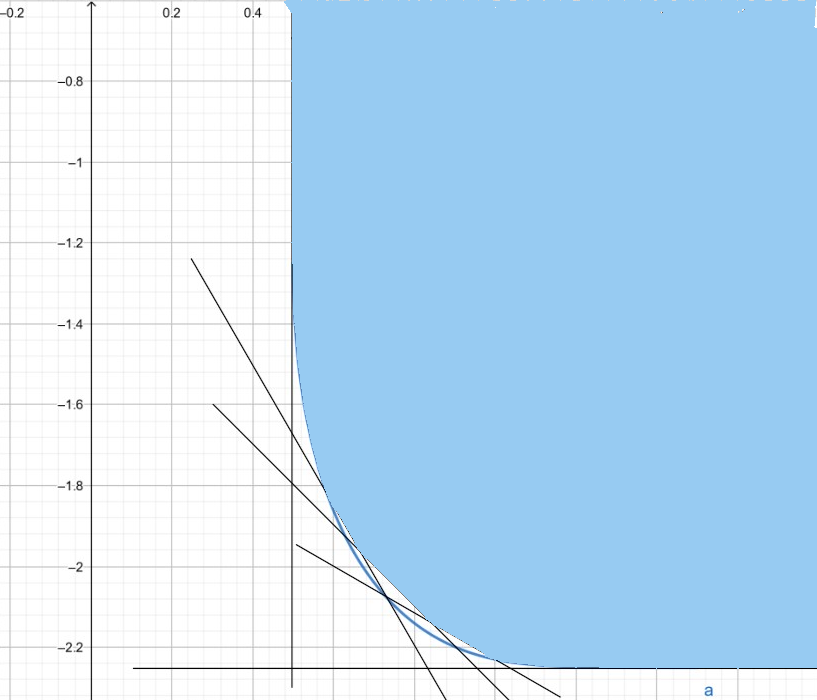}
  \end{center}
  It is possible to see that the hypothesis \eqref{hyp_U} holds and $U(t,x)$ is a 
  solution of the Hamilton-Jacobi equation.
\end{example}

In the following example hypothesis \eqref{hyp_U} does not hold.

\begin{example}
  Like before, we choose $\RR^n=\RR^d=\RR^2$ and $C=\RR^2_+$.  The 
  Hamiltonian and the initial condition are 
  $$
  \Ha_0(p)=\left(\begin{array}{cc}
    \frac12\| p\|^2 \\ \frac12\| p+p_0\|^2 \end{array}\right)+C, \quad 
  \Ha(p,\zeta)=\Ha_0(p)+C,\quad
  U_0(x)=\left(\begin{array}{cc}
    \frac12\| x\|^2 \\ \frac12\| x\|^2 \end{array}\right).
  $$
  For $\zeta={\zeta_1\choose1-\zeta_1}\in B^+(\hat z)$, where $\hat z$ is the same 
  as in the previous example, we have
  $$
  \begin{aligned}
  \Ha_\zeta(p) &= \frac{\zeta_1}{2}\| p\|^2 +\frac{1-\zeta_1}{2}\| p+p_0\|^2, \\
  D\Ha_\zeta(p) &= p+(1-\zeta_1) p_0, \\
  U_{0,\zeta}(x) &= \frac{1}{2}\| x\|^2, \\
  DU_{0,\zeta}(x) &= x.
  \end{aligned}
  $$
  The solutions \eqref{XVP} are
  $$
  \left\{ \begin{array}{l}
  X_\zeta(t,x) = (1+t)x+t(1-\zeta_1) p_0 \\
  V_\zeta(t,x) = \frac{1}{2}(1+t)\| x\|^2-\frac{1-\zeta_1}{2}t\| p_0\|^2 \\
  P_\zeta(t,x) = x
  \end{array} \right.
  $$
  We find that $Z_\zeta(t,x)$ is well defined for any nonnegative $t$ and
  $$
  Z_\zeta(t,x)=\frac{x-t(1-\zeta_1) p_0}{1+t}.
  $$
  The solutions $u_\zeta(t,x)$ of \eqref{Ham_Jac_eq_scal} are global for any 
  $\zeta$ and any $x\in\RR^n$.  Correspondingly, we obtain
  $$
  U_\zeta(t,x) =\left(\begin{array}{c}
  \frac{1}{2(1+t)}\| x-t(1-\zeta_1) p_0\|^2 \\
  \frac{1}{2(1+t)}\| x-t(1-\zeta_1) p_0\|^2-\frac{1}{2}t\| p_0\|^2
  \end{array}\right)+H^+(\zeta).
  $$
  To check the property \eqref{hyp_U}, the curve that describes $U_\zeta(t,x)$ is 
  plotted for $\zeta\in B^+(\hat z)$, $x=p_0=\binom{1}{0}$ and $t=1$:\\
  \begin{center}
  \includegraphics[width=13cm]{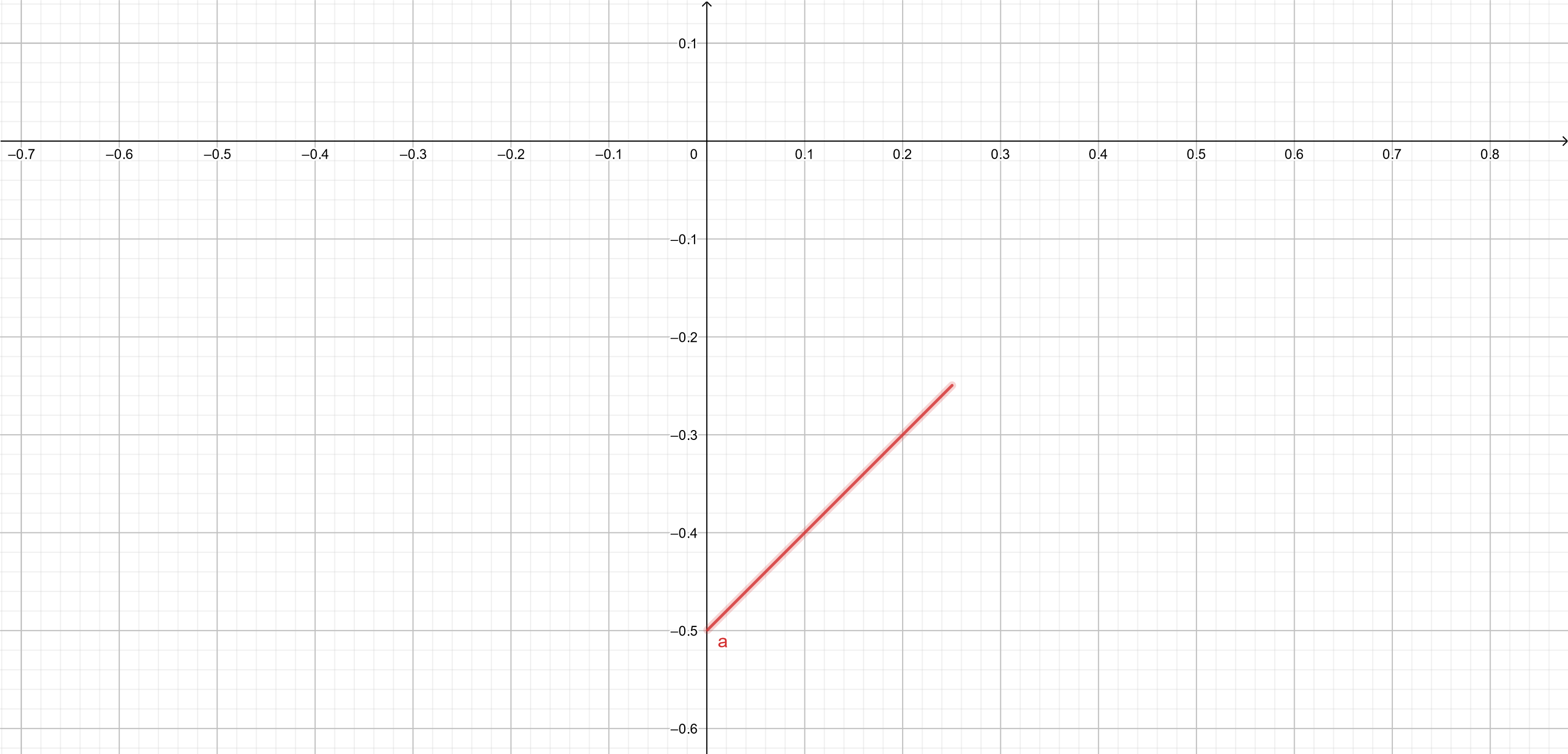}\\
  \end{center}
  In the following figures the half-spaces corresponding to $\zeta={1\choose0}$, 
  ${1/2\choose1/2}$ and ${0\choose1}$ are drawn.  It is possible to observe that 
  in the second figure the half-space corresponding to $\zeta={1/2\choose1/2}$ is 
  not on the border of the intersection of the other two half-spaces, so hypothesis 
  \eqref{hyp_U} is not fulfilled.
  \begin{center}
  \includegraphics[width=6.5cm]{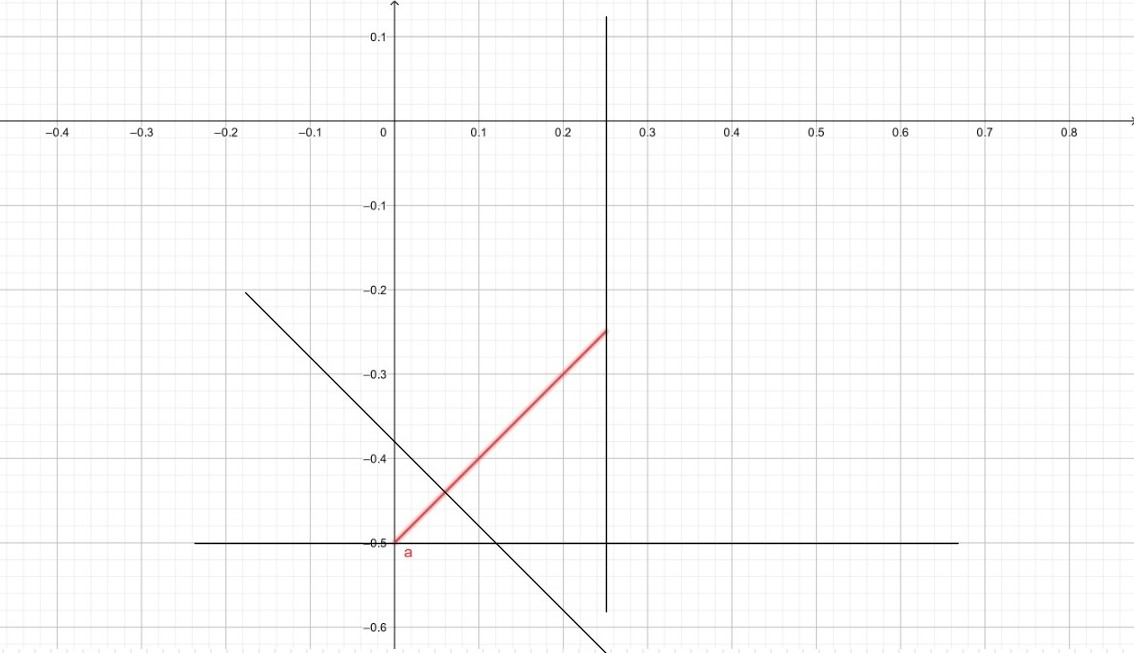}\quad
  \includegraphics[width=6.5cm]{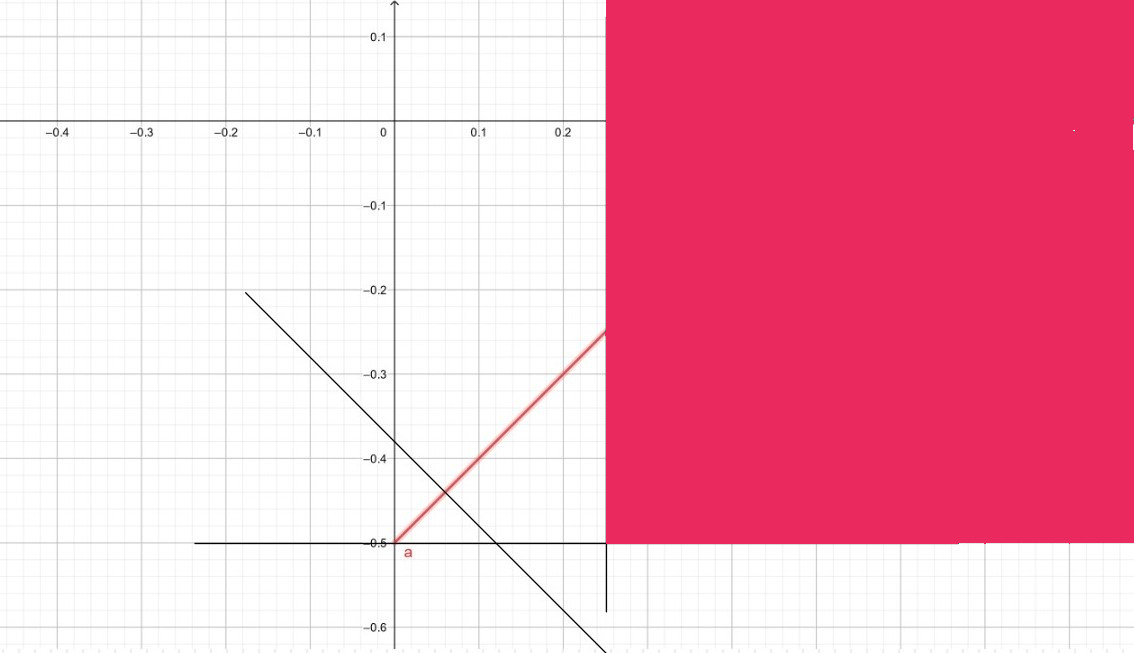}
  \end{center}
\end{example}

\section{A scalarization solution for a multiobjective calculus of variations problem}

Let us consider the continuous lower bounded functions
$$
L:\RR^n\to\RR^d\, ,\qquad U_0:\RR^n\to\RR^d
$$
where $L$ is the running cost or Lagrangian and $U_0$ is the initial cost.

For any $(t,x)\in[0,+\infty)\times\RR^n$, define the set of admissible arcs:
$$
Y(t,x)=\{ y\in W^{1,1}([0,t],\RR^n)\ \vert\ y(t)=x\}\, .
$$
In \cite{hv} the problem of `minimizing'  the cost functional $J_t:W^{1,1}([0,t],
\RR^n)\to\RR^d$
$$
J_t[y]=\int_0^t L(s,y(s),\dot y(s))\ ds + U_0(y(0))
$$
with respect to $y\in Y(t,x)$ was considered.

In order to precise the meaning of the previous minimization, we consider 
the functions:
\begin{eqnarray*}
\overline L&:&\RR^n\to\G(\RR^d,C)\\
\overline J_t &:&W^{1,1}([0,t],\RR^n)\to\G(\RR^d,C)
\end{eqnarray*}
defined by the inf-extension $\overline L(s,y,z)=L(s,y,z)+C$ and
$$
\overline J_t[y]=\int_0^t \overline L(s,y(s),\dot y(s))\ ds + U_0(y(0)),
$$
where the integral is in the Aumann sense (see \cite{Aumann} or \cite{Aubin_Frank}).

Now the problem can be written:
\begin{equation}
\label{pb}
\mbox{minimize }\overline J_t[y]\mbox{ over all arcs }y\in Y(t,x)\, .
\end{equation}
Since the functional $\overline J_t$ maps into the complete lattice 
$\G(\RR^d,C)$, the value function is well defined:
\begin{equation}
\label{value_f}
U(t,x)=\inf_{y\in Y(t,x)} \overline J_t[y].
\end{equation}

Let $\overline L:\RR^n\to\G(\RR^d,C)$ be a convex function.  For any $\zeta
\in B^+(\hat z)$ let $L_\zeta(p)=L(p)\cdot\zeta$ be such that
\begin{equation}
\label{coerc}
\lim_{\|p\|\to+\infty}\frac{L_\zeta(p)}{\| p\|}=+\infty
\end{equation}
and let $U_{0,\zeta}(x)=U_0(x)\cdot\zeta$ be Lipschitz continuous.

For any $\zeta\in B^+(\hat z)$ there exists $w_\zeta$ such that
\begin{equation}
\label{w_zeta}
\inf_{w\in\RR^n} \left[ tL_\zeta\left(\frac{x-w}{t}\right)+U_{0,\zeta}(w)\right]
  =\left[ tL_\zeta\left(\frac{x-w_\zeta}{t}\right)+U_{0,\zeta}(w_\zeta)\right].
\end{equation}
The element $w_\zeta$ is a $\zeta$-minimizer.

The value function $U(t,x)$ was proved to be obtained as an infimum over $\RR^n$ 
through the Hopf-Lax formula:
$$
U(t,x)=\inf_{w\in\RR^n} \left[ t\overline L\left(\frac{x-w}{t}
  \right)+U_{0}(w)\right].
$$

We prove now that it is sufficient to take the infimum over a smaller set, instead of all 
$\RR^n$.  More precisely, one can consider only the set of the $\zeta$-minimizers, 
for $\zeta\in B^+(\hat z)$.

\begin{theorem}
  Let $L:\RR^n\to\RR^d$ and $U_0:\RR^n\to\RR^d$ be continuous functions, 
  $\overline L:\RR^n\to\G(\RR^d,C)$ be a convex function and satisfy \eqref{coerc}.  
  Then the value function $U$ with values in $\G(\RR^d,C)$ is given by the formula
  $$
  U(t,x)=\inf_{\zeta\in B^+(\hat z)} \left[ t\overline L\left(\frac{x-w_\zeta}{t}
    \right)+U_{0}(w_\zeta)\right].
  $$
\end{theorem}

\begin{proof}
  If we denote
  $$
  V(t,x)=\inf_{\zeta\in B^+(\hat z)} \left[ t\overline L\left(\frac{x-w_\zeta}{t}
    \right)+U_{0}(w_\zeta)\right],
  $$
  it is immediate that $V(t,x)\subseteq U(t,x)$.  Let $z_0\in U(t,x)\bs V(t,x)$.  Since 
  $V(t,x)$ is closed and convex, by the separation theorem there exists a nonzero 
  $\xi\in\RR^d$ such that
  \begin{equation}
  \label{ineq_sep}
  \xi\cdot z_0<K<\xi\cdot z
  \end{equation}
  for any $z\in V(t,x)$.  We want to show that $\xi\in C^+$.  In fact, if not 
  there exists $c\in C$ with $\xi\cdot c<0$.  Now, if $z\in V(t,x)$, also $z+\lambda c
  \in V(t,x)$ for any $\lambda\geq 0$.  The following limit holds
  $$
  \lim_{\lambda\to+\infty} \xi\cdot(z+\lambda c)=-\infty,
  $$
  but this contradicts inequality \eqref{ineq_sep}.  It is always possible to consider 
  $\xi$ in $B^+(\hat z)$.  Now \eqref{ineq_sep} implies that
  $$
  \xi\cdot z_0<tL_\xi\left(\frac{x-w_\xi}{t}\right)+U_{0,\xi}(w_\xi)
  $$
  and this is not possible.
\end{proof}

The previous theorem clarifies also that the set of all the linear arcs
$$
y_\zeta(s)=w_\zeta+\frac{s}{t}(x-w_\zeta)
$$
for $\zeta\in B^+(\hat z)$ forms an infimizer for problem \eqref{pb} and 
more precisely a scalarization solution:

\begin{corollary}
  The set
  \begin{equation}
  \label{scalinf}
  M=\{ y_\zeta\in Y(t,x)\mid\zeta\in B^+(\hat z)\}
  \end{equation}
  is a scalarization solution for problem \eqref{pb}.
\end{corollary}

\section{Properties of the set-valued Fenchel conjugate}

In the following lemma and theorem some properties of the Fenchel conjugate 
are stated. In the lemma the link between the set-valued and the scalarized 
Fenchel conjugate is studied. 

\begin{lemma}
  For any $\zeta\in B^+(\hat z)$ the following equalities hold:
  \begin{eqnarray}
  \overline L^*(p,\zeta) & = & S_{(1,\zeta)}(L_\zeta^*(p))\label{L*set} \\
  \inf_{z\in\overline L^*(p,\zeta)}\zeta\cdot z & = & L_\zeta^*(p). \label{L*scal}
  \end{eqnarray}
\end{lemma}

\begin{proof}
  The Fenchel conjugate $\overline L^*(p,\zeta)$ is defined as the supremum over 
  $\RR^n$ of $S_{(p,\zeta)}(x)-_\zeta \overline L(x)$.  Each of the half-spaces can 
  be written as
  $$
  S_{(p,\zeta)}(x)-_\zeta \overline L(x)=(p\cdot x-L_\zeta(x))\hat z
    +H^+(\zeta).
  $$
  Since the half-spaces are parallel, we have
  $$
  \begin{aligned}
  \sup_{x\in\RR^n}\left[S_{(p,\zeta)}(x)-_\zeta \overline L(x)\right] &=
    \bigcap_{x\in\RR^n}\left[S_{(p,\zeta)}(x)-_\zeta \overline L(x)\right] \\
    &=\left[\sup_{x\in\RR^n}(p\cdot x-L_\zeta(x))\right]\hat z+H^+(\zeta) \\
    &=L_\zeta^*(p)\hat z+H^+(\zeta).
  \end{aligned}
  $$
  This proves \eqref{L*set} and \eqref{L*scal}.
\end{proof}

In the assumption that $L_\zeta$ is $C^2$, coercive (hypothesis \eqref{coerc}) 
and strictly convex, some well-known properties of the Fenchel conjugate (see 
for example \cite{Cannarsa_Sinestrari}) hold:
\begin{equation}
\label{Fenchel_class}
\begin{aligned}
DL^*_\zeta(p_0) &=(DL_\zeta)^{-1}(p_0), \\
D^2L^*_\zeta(p_0) &= \left[ D^2L_\zeta(DL^*_\zeta(p_0))\right]^{-1}, \\
L^*_\zeta(p_0) &= p_0\cdot DL^*_\zeta(p_0)-L_\zeta(DL^*_\zeta(p_0)).
\end{aligned}
\end{equation}
In the next theorem the previous properties are extended to the set-valued case.

\begin{theorem}
  Given $\zeta\in B^+(\hat z)$, suppose that $L$ is of 
  class $C^2$, satisfies \eqref{coerc}  and $L_\zeta(p)$ is strictly 
  convex.   Then $\overline L^*(p,\zeta)$ is twice differentiable in $p$ with respect 
  to $\zeta$ and
  \begin{eqnarray}
  D_{\zeta,p} \overline L^*(p_0,\zeta)(p) & = & S_{((DL_\zeta)^{-1}(p_0),\zeta)}
      (p), \label{DL*} \\
  D^2_{\zeta,p}\overline L^*(p_0,\zeta)(p_1,p_2) & = & S_{\left(p_1^T[D^2L_\zeta
      (DL^*_\zeta(p_0))]^{-1},\zeta\right)}(p_2) \label{D2L*} \\
  \overline L^*(p_0,\zeta) & = & S_{\left(DL^*_\zeta(p_0),\zeta\right)}(p_0)
      -_\zeta\overline L\left(DL^*_\zeta(p_0)\right),
  \end{eqnarray}
  where $p_1^T$ is the transpose of the vector $p_1$.
\end{theorem}

\begin{proof}
  In order to calculate the first derivative of $\overline L^*(\cdot,\zeta)$ at $p_0
  \in\RR^n$ in the direction $p\in\RR^n$ with respect to $\zeta$, we must study 
  the limit
  $$
  \lim_{h\to 0^+} \frac1h\left[ \overline L^*(p_0+hp,\zeta)-_\zeta\overline L^*
    (p_0,\zeta)\right].
  $$
  Using the previous lemma, we obtain
  $$
  \begin{aligned}
  \frac1h &\left[ \overline L^*(p_0+hp,\zeta)-_\zeta\overline L^*
      (p_0,\zeta)\right] \\
    &=\left\{ z\in\RR^d\mid\zeta\cdot z\geq\frac1h\left[
      \inf_{z_1\in\overline L^*(p_0+hp,\zeta)}\zeta\cdot z_1-\inf_{z_2\in
      \overline L^*(p_0,\zeta)}\zeta\cdot z_2\right]\right\} \\
    &=\left\{ z\in\RR^d\mid\zeta\cdot z\geq\frac1h\left[L^*_\zeta(p_0+hp)-
      L^*_\zeta(p_0)\right]\right\}
  \end{aligned}
  $$
  and it is possible to calculate the limit
  $$
  D_{\zeta,p} \overline L^*(p_0,\zeta)(p)=\left\{ z\in\RR^d\mid\zeta\cdot z\geq 
    DL^*_\zeta(p_0)\cdot p\right\}.
  $$
  The first equation in \eqref{Fenchel_class} completes the proof of \eqref{DL*}.
  
  In order to study the second derivative, we calculate
  $$
  \begin{aligned}
  \frac1h &\left[D_{\zeta,p} \overline L^*(p_0+hp_2,\zeta)(p_1)
      -_\zeta D_{\zeta,p} \overline L^*(p_0,\zeta)(p_1)\right] \\
    &= \frac1h\left[ S_{((DL_\zeta)^{-1}(p_0+hp_2),\zeta)}(p_1)
      -_\zeta S_{((DL_\zeta)^{-1}(p_0),\zeta)}(p_1)\right].
  \end{aligned}
  $$
  Using the first equation in \eqref{Fenchel_class}, we obtain
  $$
  \begin{aligned}
  \frac1h&\left[ S_{((DL_\zeta)^{-1}(p_0+hp_2),\zeta)}(p_1)
      -_\zeta S_{((DL_\zeta)^{-1}(p_0),\zeta)}(p_1)\right] \\
    &=\frac1h \left[ S_{(DL^*_\zeta(p_0+hp_2),\zeta)}(p_1)
      -_\zeta S_{(DL^*_\zeta(p_0),\zeta)}(p_1)\right] \\
    &= \left\{ z\in\RR^d\mid \zeta\cdot z\geq\frac1h\left[DL^*_\zeta(p_0+hp_2)
      -DL^*_\zeta(p_0)\right]\cdot p_1\right\}
  \end{aligned}
  $$
  and, taking the limit,
  $$
  D^2_{\zeta,p}\overline L(p_0)(p_1,p_2)=\left\{ z\in\RR^d\mid \zeta\cdot z\geq 
    p_1^TD^2L^*_\zeta(p_0)p_2\right\}
  $$
  and this with the second equation in \eqref{Fenchel_class} completes the proof 
  of \eqref{D2L*}.
  
  By \eqref{L*set}, we have
  $$
  \overline L^*(p_0,\zeta)=S_{(1,\zeta)}(L^*_\zeta(p_0)).
  $$
  Using the third property in \eqref{Fenchel_class}, we may write
  $$
  \begin{aligned}
  \overline L^*(p_0,\zeta) &=S_{(1,\zeta)}(p_0\cdot DL^*_\zeta(p_0)-
      L_\zeta(DL^*_\zeta(p_0))) \\
    &= S_{(DL^*_\zeta(p_0),\zeta)}(p_0) -_\zeta\overline L(DL^*_\zeta(p_0)).
  \end{aligned}
  $$
\end{proof}

\begin{remark}
  If we apply the Fenchel conjugate twice
  $$
  \overline L^{**}(p,\zeta)=\sup_{q\in\RR^n}[S_{(p,\zeta)}(q)-_\zeta \overline L^*
    (q,\zeta)]
  $$
  for $\zeta\in B^+(\hat z)$ and if $\overline L$ is convex, it is easy 
  to see that
  $$
  \bigcap_{\zeta\in B^+(\hat z)} \overline L^{**}(p,\zeta)=\overline L(p).
  $$
  See for example \cite{SetOptSurvey} for a generalization of the Fenchel-Moreau 
  theorem.
\end{remark}

\section{The scalarization solution and the characteristic curves}

Let $L:\RR^n\to\RR^d$ and $U_0:\RR^n\to\RR^d$ be continuous functions, 
$\overline L:\RR^n\to\G(\RR^d,C)$ be a convex function and satisfy \eqref{coerc}. 
The value function \eqref{value_f} of the minimization problem \eqref{pb} has 
the property \eqref{u_zeta}, as it is stated in the following lemma.

\begin{lemma}
  The value function \eqref{value_f} of the minimization problem \eqref{pb} 
  is such that for any $\zeta\in B^+(\hat z)$
  $$
  U(t,x)+H^+(\zeta)=u_\zeta(t,x)\zeta+H^+(\zeta),
  $$
  where
  $$
  u_\zeta(t,x)=tL_\zeta\left(\frac{x-w_\zeta}{t}\right)+U_{0,\zeta}(w_\zeta),
  $$
  for $w_\zeta$ as in \eqref{w_zeta}.
\end{lemma}

\begin{proof}
  Recalling the Hopf-Lax formula, one has that
  $$
  \begin{aligned}
  u_\zeta(t,x) &=\inf_{z\in U(t,x)}\zeta\cdot z=\inf_{w\in\RR^n} \left[tL_\zeta\left(
      \frac{x-w}{t}\right)+U_{0,\zeta}(w)\right] \\
    &=tL_\zeta\left(\frac{x-w_\zeta}{t}\right)+U_{0,\zeta}(w_\zeta).
  \end{aligned}
  $$
\end{proof}

The value function was proved in \cite{hv} to satisfy a Hamilton-Jacobi equation 
and we report here the result:

\begin{theorem}
  Let $(t,x)\in [0,+\infty)\times\RR^n$, $\zeta\in C^+$, $\|\zeta\|=1$.  Let 
  $\overline L$ be convex, \eqref{coerc} be satisfied, $U_{0,\zeta}$ be 
  Lipschitz on $\RR^n$ and $L$, $U_0$ be of class $C^2$.  If $w_\zeta$ 
  is as in \eqref{w_zeta}, let the sum of the hessian matrices
  $$
  \frac1t H_{L_\zeta}\left(\frac{x-w_\zeta}{t}\right)+H_{U_{0,\zeta}}(w_\zeta)
  $$
  be non-singular.
  
  Then the value function $U(t,x)$ is a solution of the Hamilton-Jacobi equation
  $$
  U_{t,\zeta}(t,x)+\overline L^*(Du_\zeta(t,x),\zeta)=H^+(\zeta)
  $$
\end{theorem}

The following proposition shows the link between the scalarization solution of the 
calculus of variations problem \eqref{pb} described in \eqref{scalinf} and the 
characteristic curves for the corresponding Hamilton-Jacobi equation.

\begin{proposition}
  For any $\zeta\in B^+(\hat z)$ suppose that $L$ is 
  $C^2(\RR^n;\RR^d)$, satisfies \eqref{coerc}  and $L_\zeta(p)$ is strictly 
  convex.   Suppose also that the functions $U_{0,\zeta}$ are Lipschitz and of 
  class $C^1$.  The elements $y_\zeta$ of the set $M$ defined in \eqref{scalinf} 
  are such that
  $$
  y_\zeta(s)=X_\zeta(s,w_\zeta),
  $$
  with $w_\zeta$ as in \eqref{w_zeta}.
\end{proposition}

This results shows that there is a scalarization solution which is formed by 
characteristic curves.

\end{document}